\documentclass[graybox]{svmult}

\usepackage{amsfonts,amsmath,amstext,latexsym,amssymb}
\usepackage{mathptmx}       
\usepackage{helvet}         
\usepackage{courier}        
\usepackage{type1cm}        
\usepackage{makeidx}         
\usepackage{graphicx}        
\usepackage{multicol}        
\usepackage[bottom]{footmisc}

\usepackage{cite}

\usepackage[inline]{enumitem}  
\setlist{labelindent=1pt,itemsep=.5em}
\setlist[itemize]{leftmargin=1cm}
\setlist[enumerate]{itemindent=0em,leftmargin=1cm}

\usepackage[draft=false,hidelinks,colorlinks=true,linkcolor=blue,citecolor=black,urlcolor=black,anchorcolor=blue,bookmarks=false,hypertexnames=true,backref=page]{hyperref}

\allowdisplaybreaks[3]

\makeindex             

\newtheorem{rem:DjinjaTumwSilv1}{Remark}[section]

\newtheorem{cor:DjinjaTumwSilv1}[thm]{Corollary}
\newtheorem{lem:DjinjaTumwSilv1}[thm]{Lemma}
\newtheorem{prop:DjinjaTumwSilv1}[thm]{Proposition}
\newtheorem{ex:DjinjaTumwSilv1}[thm]{Example}

\DeclareMathOperator{\esssup}{ess\;sup}
\DeclareMathOperator{\supp}{\rm supp\,}
\smartqed

\begin{document}

\title*{Representations of polynomial covariance type commutation relations by piecewise function multiplication and composition operators}
\titlerunning{Representations by piecewise function multiplication and composition operators}
\author{Domingos Djinja, Sergei Silvestrov, Alex Behakanira Tumwesigye}
\authorrunning{D. Djinja, S. Silvestrov, A. B. Tumwesigye}
\institute{Domingos Djinja \at Department of Mathematics and Informatics, Faculty of Sciences, Eduardo Mondlane University, Box 257, Maputo, Mozambique.
Division of Mathematics and Physics, Division of Mathematics and Physics, School of Education, Culture and Communication, M\"{a}lardalens University, M\"{a}lardalens University, Box 883, V\"{a}ster\.{a}s, Sweden. \\
\email{domingos.djindja@uem.ac.mz, domingos.celso.djinja@mdu.se}
\and Sergei Silvestrov \at Division of Mathematics and Physics, School of Education, Culture and Communication, M\"{a}lardalens University, Box 883, V\"{a}ster\.{a}s, Sweden.  \\ \email{sergei.silvestrov@mdu.se}
\and Alex Behakanira Tumwesigye \at Department of Mathematics, College of Natural Sciences, Makerere University, Box 7062, Kampala, Uganda. \\ \email{alexbt@cns.mak.ac.ug}}
%
%


\maketitle
\label{Chap:DjinjaSilvestrovTumwesigye2repscovarrelLp}

\abstract*{Representations  of polynomial covariance type commutation relations
are constructed on  Banach spaces $L_p$ and $C[\alpha, \beta],\ \alpha,\beta\in \mathbb{R}$.
Representations involve operators with piecewise functions,
 multiplication operators and inner superposition operators.
 \keywords{piecewise function multiplication operator, composition operator, covariance commutation relations}\\
{\bf MSC2020 Classification:}  47B33, 47L80, 47L65
}

\abstract{Representations  of polynomial covariance type commutation relations
are constructed on  Banach spaces $L_p$ and $C[\alpha, \beta]\ \alpha,\beta\in \mathbb{R}$.
Representations involve operators with piecewise functions,
multiplication operators and inner superposition operators.
 \keywords{piecewise function multiplication operator, composition operator, covariance commutation relations}\\
{\bf MSC2020 Classification:} 47B33, 47L80, 47L65
}






\section{Introduction}         \index{commutation relations}
In many areas of applications, there can be found relations of the form
\begin{gather} \label{covrelation1}
  ST=F(TS)
\end{gather}
where $S, T$ are elements of
an associative algebra and   $F:\mathbb{R}\to\mathbb{R}$ is a function satisfying certain conditions.
For example, if $F(z)=z$, then $S$ and $T$ commute. If $F(z)=-z$, then $S$ and $T$ anti-commute.
If $F(z)=\delta_0 +\delta_1 z$,
then $S$ and $T$ satisfy the relation
\begin{gather}\label{deformedHeisembergCR}
ST-\delta_1 TS=\delta_0 I ,
\end{gather}
where $\delta_0,\delta_1$ are constants and $I$ is the identity element. This relation is known as
the deformed Heisenberg commutation relation \cite{HellstromSSDST2}.
If $\delta_0=\delta_1=1$, then \eqref{deformedHeisembergCR} reduces to the canonical Heisenberg
commutation relation
$
ST-TS=I
$
which is important in differential and integral calculus and quantum physics.
If $\delta_0=0$ then \eqref{deformedHeisembergCR} reduces to the quantum plane relation
$
ST=\delta_1 TS.
$
Elements of an algebra that satisfy a commutation relation are called a representation of this commutation relation in that algebra. Representations of covariance commutation relations \eqref{covrelation1} by linear operators are important for study of actions and induced representations of groups and semigroups, crossed product operator algebras, dynamical systems, harmonic analysis, wavelets and fractals analysis and applications in physics and engineering
\cite{BratJorgIFSAMSmemo99,BratJorgbook,DaubechiesIngridDST2,JorgWavSignFracbook,JorgOpRepTh88,JorMoore84,MACbook1,MACbook2,MACbook3,OstSambookDST2,Pedbook79,SamoilenkobookDST2}.
The structure of representations for the relations of the form \eqref{covrelation1}, \eqref{covrelationNormal}, \eqref{covrelationReverse},
by bounded and unbounded self-adjoint operators, normal operators, unitary operators, partial isometries
and other linear operators with special involution conditions on a Hilbert space, have been considered in   \cite{BratEvansJorg2000,CarlsenSilvExpoMath07,CarlsenSilvAAM09,CarlsenSilvProcEAS10,DutkayJorg3,DJS12JFASilv,DLS09,DutSilvProcAMS,
DutSilvSV,JSvT12a,JSvT12b,Mansour16,JMusondaPhdth18,JMusonda19,Musonda20,Nazaikinskii96,OstSambookDST2,PerssonSilvestrov031,
PerssonSilvestrov032,PersSilv:CommutRelDinSystDST2,RST16,RSST16,SamoilenkobookDST2,SaV8894,SilPhD95,STomdynsystype1,
SilWallin96,SvSJ07a,SvSJ07b,SvSJ07c,SvT09,Tomiyama87,Tomiama:SeoulLN1992,Tomiama:SeoulLN2part2000,AlexThesis2018DST2,TumwRiSilv:ComCrPrAlgPieccnstfnctreallineSPAS19v2,VaislebSa90}
using reordering formulas for functions of the algebra elements and operators satisfying covariance commutation relation,
functional calculus and spectral representation of operators and interplay with dynamical systems generated
by iteration of the maps involved in the commutation relations.

In case of $*$-algebras, when $S=X$ and $T=X^*$ where $X$ is an element in the algebra, the relation \eqref{covrelation1} reduces to $ XX^*=F(X^*X)$.
This relation often can be transformed to relations of the form
\begin{align}\label{covrelationNormal}
& AB=BF(A),\\
\label{covrelationReverse}
& BA=F(A)B
\end{align}
for some other elements $A,B$ of the $*$-algebra obtained from $X$ and $X^*$ using some transformations or factorizations in an appropriate functional calculus  (see for example \cite{OstSambookDST2,PersSilv:CommutRelDinSystDST2,SamoilenkobookDST2,SaV8894,VaislebSa90,AlexThesis2018DST2} and references cited their).
A description of the structure of representations of relation \eqref{covrelationNormal} by bounded and unbounded self-adjoint linear operators on Hilbert space by using spectral representation \cite{AkkLnearOperatorsDST2} of such operators is given in \cite{SamoilenkobookDST2}, where also more general families of commuting self-adjoint operators satisfying relation \eqref{covrelationNormal} with other operators on Hilbert spaces are considered using spectral theory and non-commutative analysis for bounded and unbounded operators on Hilbert spaces.

In this paper we construct representations of relation \eqref{covrelationNormal} and \eqref{covrelationReverse} by linear operators acting on Banach spaces $L_p$ and $C[\alpha,\beta]$ for $\alpha,\beta\in \mathbb{R}$, and $F$ is a polynomial. Such representations
can also be viewed as solutions for operator equations $AX=XF(A)$, when $A$ is specified or
$XB=BF(X)$ when $B$ is specified. We consider representations of \eqref{covrelationNormal} involving linear operators with piecewise function multiplication operators, multiplication operators and inner superposition operators. We derive conditions on the parameters or functional coefficients of operators  so that they satisfy
\eqref{covrelationNormal}
for a polynomial $F$. In contrast to \cite{OstSambookDST2,SamoilenkobookDST2,SaV8894,VaislebSa90} devoted to involutive representations of covariance type relations by operators on Hilbert spaces using spectral theory of operators on Hilbert spaces,
aiming at direct construction of various classes of representations of covariance type relations in specific important classes of operators on Banach spaces more general than Hilbert spaces without imposing any involution conditions and not using classical spectral theory of operators.

This paper is organized in five sections, after the introduction,  we present in Section  \ref{SecPreNotDST2} preliminaries, notations and basic definitions.
In Section \ref{SecOpPiecewiseFunctDST2}, we present representations involving  piecewise function multiplication operators acting on $L_p$ for $1<p<\infty$. In Section \ref{SecInnerSuperOpDST2}, we construct representations involving inner superposition operators. These operators are important in wavelets analysis for instance.
In Section \ref{SecMultiplOpDST2} we construct representations by multiplication operators acting on $L_p$ for $1<p<\infty$  and the space of continuous functions.

\section{Preliminaries and notations}\label{SecPreNotDST2}
We use the following basic standard definitions and notations
(see for example \cite{AdamsGDST2,BrezisFASobolevSpacesDST2,ConwayFunctionalAnalysisDST2,DudleyDST2,FollandRADST2,KantarovitchAkilovFunctionalAnalysisDST2,KolmogorovDST2,RudinRCADST2}).
Let $S\subseteq \mathbb{R}$, ($\mathbb{R}$ is the set of real numbers), be a Lebesgue measurable set and let $(S,\Sigma, \tilde{m})$ be a $\sigma$-finite measure space, that is, $S$ is a nonempty set, $\Sigma$ is a $\sigma-$algebra with subsets of $S$,
where $S$ can be cover with at most  countably many disjoint sets $E_1,E_2,E_3,\ldots$ such that $ E_i\in \Sigma, \,
\tilde{m}(E_i)<\infty$, $i=1,2,\ldots$  and $\tilde{m}$ is the Lebesgue measure.
For $1\leqslant p<\infty,$ we denote by $L_p(S)$, the set of all classes of equivalent measurable functions $f:S\to \mathbb{R}$ such that
$\int\limits_{S} |f(t)|^p dt < \infty.$
This is a Banach space with norm
$\| f\|_p= \left( \int\limits_{S} |f(t)|^p dt \right)^{\frac{1}{p}}.$
We denote by $L_\infty(S)$ the set of all classes of equivalent measurable functions $f:S\to \mathbb{R}$ such that there exists a constant $\lambda>0$ for which
$|f(t)|\leq \lambda$
for almost every $t$. This is a Banach space with norm
$ \|f\|_{\infty}=\mathop{\esssup}_{t\in S} |f(t)|.$
We denote by $C[\alpha,\beta]$ the set of all continuous functions $f:[\alpha,\beta]\to \mathbb{R}$. This is a Banach
space with norm
$\|f\|=\max_{t\in [\alpha,\beta]} |f(t)|.$

\section{Representations by operators involving piecewise functions}\label{SecOpPiecewiseFunctDST2}
When operators are given in abstract form,  we use the notation $A:L_p(\mathbb{R})\to L_p(\mathbb{R})$ meaning that operator $A$  is well defined
   from $L_p(\mathbb{R})$ to $L_p(\mathbb{R})$ without discussing sufficient conditions for that to be satisfied.

\begin{prop:DjinjaTumwSilv1}\label{PropPiecewiseFunct}
Let $(\mathbb{R},\Sigma,\tilde{m})$ be the standard Lebesgue measure space in the real line.
Let $A:L_p(\mathbb{R})\to L_p(\mathbb{R})$,  $B:L_p(\mathbb{R})\to L_p(\mathbb{R})$, $1\le p\le\infty$ be defined as follows
\begin{gather}\label{DefinitionOfOperatorsAandBPiecewiseFunct}
   (Ax)(t)=\sum_{i=1}^{l} \alpha_i a(t)I_{G_i}(t) x(t),\qquad
   (Bx)(t)=\sum_{i=1}^{m} \beta_i b(t)I_{H_i}(t) x(t),
\end{gather}
where $a,b:\mathbb{R}\to\mathbb{R} $ are measurable functions, $G_i\in \Sigma$ for $i=1,\ldots, l$, $\tilde{m}(G_i\cap G_j)=0$ if $i\not=j$,
$H_i\in \Sigma$ for $i=1,\ldots, m$, $\tilde{m}(H_i\cap H_j)=0$ if $i\not=j$,
and $I_{G_k}$ is the indicator function
of the set $G_k$. Then,  for a polynomial $F_1(z)=\delta_0+F(z)$
with $F(z)=\sum\limits_{j=1}^{n}\delta_j z^j$ and
$\delta_j\in \mathbb{R}$ for $j=0,\dots,n$,
the commutation relation $AB=BF_1(A)$ is satisfied if and only if, for almost every $t$,
\begin{equation*}
  \sum_{i=1}^l \sum_{j=1}^{m} \alpha_i \beta_j a(t)b(t)I_{G_i\cap H_j}(t) =\sum_{j=1}^m \delta_0\beta_jb(t)I_{H_j}(t)+\sum_{i=1}^l\sum_{j=1}^{m} \beta_j b(t)I_{G_i\cap H_j}(t) F(a(t)\alpha_i).
\end{equation*}
\end{prop:DjinjaTumwSilv1}

\begin{proof}
  Consider a monomial $M(z)=\delta z^r$, where $r$ is a positive integer, and $\delta\in\mathbb{R}$.
  We compute $AB$, $\delta A^r$ and $\delta B A^r$ as follows:
  \begin{gather*}
    (ABx)(t)=\sum_{i=1}^{l} \alpha_i a(t)I_{G_i}(t)\left(\sum_{j=1}^{m} \beta_j b(t)I_{H_j}(t)x(t)\right)\\
    =
    \sum_{i=1}^l\sum_{j=1}^{m} \alpha_i \beta_j a(t)b(t)I_{G_i\cap H_j}(t) x(t),
    \\
    (A^2x)(t)=\sum_{i=1}^l\sum_{j=1}^{l}\alpha_i a(t)\alpha_j a(t)I_{G_i\cap G_j}(t)x(t)=\sum_{i=1}^{l} \alpha^2_ia(t)^2 I_{G_i}(t)x(t)
    \end{gather*}
   for almost all $t$.
   In general,  for almost all $t\in\mathbb{R}$,
   \begin{eqnarray*}
    (\delta A^r x)(t)&=&\sum_{i=1}^{l} \delta\alpha^r_i a(t)^r I_{G_i}(t)x(t),
    \\
    ( B \delta A^r x)(t)&=&\sum_{j=1}^m\sum_{i=1}^{l} \beta_j b(t) \delta \alpha^r_i a(t)^r I_{G_i}(t)^r I_{H_j}(t) x(t)
    \\
    &=&
    \sum_{j=1}^m\sum_{i=1}^{l} \beta_j b(t) \delta \alpha^r_i a(t)^r  I_{G_i\cap H_j}(t) x(t).
  \end{eqnarray*}
   Suppose now that $F_1(z)=\delta_0+F(z)$, where
    $F(z)=\sum\limits_{j=1}^{n}\delta_j z^j$ with
$\delta_j\in \mathbb{R}$, $j=0,\dots,n$.
   Then, for almost every $t\in\mathbb{R}$ we have
   \begin{eqnarray*}
     (F_1(A)x)(t)&=&\delta_0x(t)+(F(A)x)(t)=\delta_0x(t)+\sum_{i=1}^{l} F(\alpha_i a(t)I_{G_i}(t))x(t)\\
     B(F_1(A)x)(t)&=&\delta_0(Bx)(t)+B(F(A)x)(t)\\
     &=&\sum_{j=1}^m \delta_0\beta_j b(t)I_{H_j}(t)x(t)+\sum_{i=1}^l\sum_{j=1}^{m} b(t) F(\alpha_i a(t)) \beta_j I_{G_i\cap H_j}(t)x(t).
   \end{eqnarray*}
Then $AB=BF_1(A)$ on $L_p(\mathbb{R})$ if and only if,  for almost every $t\in\mathbb{R}$,
\begin{align*}
& \sum_{i=1}^l\sum_{j=1}^{m} \alpha_i a(t)b(t) \beta_j  I_{G_i\cap H_j}(t)\\
& =\sum_{j=1}^m\delta_0 b(t)\beta_j I_{H_j}(t)
+ \sum_{i=1}^l\sum_{j=1}^{m} F(\alpha_i  a(t))b(t) \beta_j I_{G_i\cap H_j}(t),
\end{align*}
since if this equality holds for almost every $t\in\mathbb{R}$, then
multiplying both sides of this equality by $x(t)$ for any $x\in L_p(\mathbb{R})$ yields $ABx=BF_1(A)x$ on $L_p(\mathbb{R})$,
and if the equality does not hold on a set $S$ of finite positive measure $0<\tilde{m}(S)<\infty $,
then $ABI_S\neq BF_1(A)I_S$ and so $AB\neq BF_1(A)$ as linear operators on $L_p(\mathbb{R})$.
\qed
\end{proof}

\begin{ex:DjinjaTumwSilv1}{\rm
Consider the operators $$A:L_p([1,3])\to L_p([1,3]),\ B:L_p([1,3])\to L_p([1,3]),\ 1\le p\le \infty$$
defined as follows, for almost every $t$,
\begin{gather*}
  (Ax)(t)=\sum_{i=1}^{3} \alpha_i I_{G_i}(t)x(t),\quad (Ax)(t)=\sum_{i=1}^{3} \beta_i I_{H_i}(t)x(t),
\end{gather*}
where
\begin{equation}\label{PartitionIntervalDefOpAandBPropOpPwFunct}
 \begin{array}{l}
 G_1=[1,\frac{3}{2}[,\ G_2=[\frac{3}{2},2],\ G_3=]2,3],\\
 H_1=[1,2], \ H_2=]2,\frac{5}{2}],\ H_3=]\frac{5}{2},3],
\end{array}
\end{equation}
$\alpha_i,$ $\beta_i\in\mathbb{R}$, $i=1,2,3$. Let $F(z)=\delta_0+\delta_3 z^3$, where $\delta_0,\delta_3\in\mathbb{R}$.
By applying Proposition \ref{propDiagOperatorLp} we have
$
  AB=B\delta_0+\delta_3 BA^3
$
if and only if
\begin{gather*}
  \sum_{i=1}^{3}\sum_{j=1}^{3} \alpha_i\beta_j I_{G_i\cap H_j}(t)=\delta_0\sum_{j=1}^{3}\beta_j I_{H_j}(t) +\sum_{i=1}^{3}\sum_{j=1}^{3} \delta_3\alpha_i^3 \beta_j I_{G_i\cap H_j}(t).
\end{gather*}
By simplifying this we have
\begin{align*}
& \beta_1\sum_{i=1}^{2}\alpha_iI_{ G_i}(t)+\alpha_3\sum_{i=2}^{3}\beta_iI_{ H_i}(t)\\
& =\sum_{i=1}^{3}\delta_0\beta_i I_{ H_i}(t) +
  \beta_1\sum_{i=1}^{2}\delta_3\alpha_i^3I_{ G_i}(t)+\alpha_3^3 \delta_3 \sum_{i=2}^{3}\beta_i I_{ H_i}(t).
\end{align*}
This is equivalent to
\begin{gather*}
  \left\{\begin{array}{c}
           \beta_1\alpha_1=\delta_0 \beta_1+\beta_1 \delta_3 \alpha_1^3 \\
           \beta_1\alpha_2=\delta_0 \beta_1+\beta_1 \delta_3 \alpha_2^3 \\
           \beta_2\alpha_3=\delta_0 \beta_2+\beta_2 \delta_3 \alpha_3^3 \\
           \beta_3\alpha_3=\delta_0 \beta_3+\beta_3 \delta_3 \alpha_3^3
         \end{array} \right.
         \Leftrightarrow
         \left\{ \begin{array}{cc}
              F(\alpha_i)\beta_1=\alpha_i\beta_1,\quad i=1,2,  \\
              F(\alpha_3)\beta_i=\alpha_3\beta_i,\quad i=2,3.
          \end{array}\right.
\end{gather*}
This is equivalent to the following:
\begin{gather*}
\left\{\begin{array}{c} \beta_1 =0, \ \text{ or } \ F(\alpha_i)=\alpha_i,\ i=1,2, \\
F(\alpha_3)=\alpha_3, \ \text{ or } \ \beta_j=0,\ j=2,3.
\end{array}\right.
\end{gather*}
That is, one of the following cases hold:
\begin{enumerate}
\item $(\beta_1,\beta_2,\beta_3)\in (\{0\}\times (\mathbb{R}\setminus\{0\})\times \mathbb{R}) \cup (\{0\}\times \mathbb{R}\times (\mathbb{R}\setminus{0})) $, \\
    $(\alpha_1,\alpha_2,\alpha_3)\in \mathbb{R}\times\mathbb{R}\times\mbox{Fix}(F)$.
\item $(\beta_1,\beta_2,\beta_3)\in ((\mathbb{R}\setminus\{0\})\times (\mathbb{R}\setminus\{0\})\times \mathbb{R}) \cup ((\mathbb{R}\setminus\{0\})\times \mathbb{R}\times (\mathbb{R}\setminus{0})) $,\\
     $(\alpha_1,\alpha_2,\alpha_3)\in \mbox{Fix}(F)\times \mbox{Fix}(F)\times \mbox{Fix}(F)$.
\item $(\beta_1,\beta_2,\beta_3)\in ((\mathbb{R}\setminus\{0\})\times \{0\}\times \{0\})$,$(\alpha_1,\alpha_2,\alpha_3)\in \mbox{Fix}(F)\times \mbox{Fix}(F)\times \mathbb{R}$,
\item $(\beta_1,\beta_2,\beta_3)\in \{0\}\times \{0\}\times \{0\}$, $(\alpha_1,\alpha_2,\alpha_3)\in \mathbb{R}^3$.
\end{enumerate}
where Fix$(F)=\{\alpha\in\mathbb{R}:\ F(\alpha)=\alpha\}$. Note that Fix$(F)\not=\emptyset$, since $deg(F(z)-z)=1$ if $\delta_3=0$, $deg(F(z)-z)=3$ if $\delta_3\not=0$,
and all coefficients of $F(z)-z=\delta_3 z^3-z+\delta_0$ are real.

}\end{ex:DjinjaTumwSilv1}

In the following corollary we consider a case where both operators $A$ and $B$ are considered on the same partition.
\begin{cor:DjinjaTumwSilv1}\label{corSamePartitionPiecwiseFunct}
Let $(\mathbb{R},\Sigma,\tilde{m})$ be the standard Lebesgue measure space in the real line.
Let $A:L_p(\mathbb{R})\to L_p(\mathbb{R})$,  $B:L_p(\mathbb{R})\to L_p(\mathbb{R})$, $1\leq p\leq\infty$ be defined as follows
\begin{gather*}
  (Ax)(t)=\sum_{i=1}^{l} \alpha_i a(t)I_{G_i}(t) x(t),\qquad
   (Bx)(t)=\sum_{i=1}^{l} \beta_i b(t)I_{G_i}(t) x(t),
\end{gather*}
where $G_i\in \Sigma$ for $i=1,\ldots, l$, $\tilde{m}( G_i \cap G_j)=0$  for $i\not=j$, $I_{G_k}$ is the indicator function
of the set $G_k$, and
$a,b:\mathbb{R}\to\mathbb{R} $ are measurable functions. Then, for a polynomial $F_1(z)=\delta_0+F(z)$,
where $F(z)=\sum\limits_{j=1}^{n}\delta_j z^j$ and
$\delta_j\in \mathbb{R}$ for $j=0,\dots,n$, the commutation relation
$AB=BF_1(A)$ is satisfied if and only if
\begin{gather*}
  \sum_{i=1}^l \alpha_i \beta_i a(t)b(t)I_{G_i}(t) =\sum_{i=1}^l b(t)F_1(a(t)\alpha_i )\beta_i I_{G_i}(t)
  \quad \mbox{for almost every $t$}.
\end{gather*}
\end{cor:DjinjaTumwSilv1}

\begin{proof}
  By applying Proposition \ref{PropPiecewiseFunct} we have $AB=BF_1(A)$ if and only if
  \begin{gather*}
    \sum_{i,j=1}^l \alpha_i a(t)b(t) \beta_j  I_{G_i\cap G_j}(t)=\sum_{j=1}^l b(t) \delta_0 \beta_j I_{G_j}(t)+ \sum_{j,i=1}^l b(t) F(\alpha_i a(t)) \beta_j I_{G_i\cap G_j}(t)
  \end{gather*}
  for almost every $t$. Since $G_i\cap G_j=\emptyset$ when $i\not=j$, the last condition becomes
    \begin{gather*}
        \sum_{i=1}^l \alpha_i a(t)b(t) \beta_i  I_{G_i}(t)=\sum_{i=1}^l b(t) \delta_0 \beta_i I_{G_i}(t)+ \sum_{i=1}^l b(t) F(\alpha_i a(t)) \beta_i I_{G_i}(t) \\
        = \sum_{i=1}^l  b(t) \beta_i I_{G_i}(t) (\delta_0+F(\alpha_i a(t)))=\sum_{i=1}^l  b(t)\beta_i I_{G_i}(t) F_1(\alpha_i a(t)). \tag*{\qed}
    \end{gather*}
\end{proof}

\begin{cor:DjinjaTumwSilv1}
Let $(\mathbb{R},\Sigma,\tilde{m})$ be the standard Lebesgue measure space in the real line.
Consider the operators $A$ and $B$ defined in \eqref{DefinitionOfOperatorsAandBPiecewiseFunct} where $\beta_j\not=0$ for all
$j\in\{1,\ldots,m\}$.
  Then,  for a constant monomial  $F_1(z)=\delta_0$, we have $AB=BF_1(A)$ if and only if
   for all $i\in \{1,\ldots, l\}$, $j\in\{1,\ldots,m\}$ the set
  \begin{align*}
  \supp (\alpha_ia(t)-\delta_0)\cap \supp b(t)\cap (G_i\cap H_j)
  \end{align*}
  has measure zero.
\end{cor:DjinjaTumwSilv1}

\begin{proof}
We have
$
  (F_1(A)x)(t)=\delta_0 x(t),
$
for almost all $t$.
Moreover, for almost every $t$
\begin{align*}
    &  (ABx)(t)=\sum_{i=1}^{l} \alpha_i a(t)I_{G_i}(t)\left(\sum_{j=1}^{m} \beta_j b(t)I_{H_j}(t)x(t)\right)\\
    &=
   \sum_{i=1}^l\sum_{j=1}^{m} \alpha_i \beta_j a(t)b(t)I_{G_i\cap H_j} x(t)\\
  & (BF_1(A)x)(t)=\sum_{j=1}^m \delta_0 \beta_j b(t) I_{H_j}(t)x(t),
\end{align*}
 Then, for almost every $t\in \mathbb{R}$ we have $AB=BF_1(A)$ if and only if
\begin{equation*}
    \sum_{i=1}^l \sum_{j=1}^{m} \alpha_i a(t) \beta_j b(t) I_{G_i\cap H_j}(t)=\sum_{j=1}^m \delta_0 \beta_j b(t) I_{H_j}(t),
  \end{equation*}
  which is equivalent to
  \begin{align}\label{cor:picewisfunctionsPolynomialConst}
    & \sum_{i=1}^l \alpha_i a(t) \beta_j b(t) I_{G_i\cap H_j}(t)= \delta_0 \beta_j b(t) I_{H_j}(t) \\
     \nonumber
    & \quad \text{for almost every}\ t\in\mathbb{R}, \ j\in \{1,\ldots,m\}
  \end{align}
Since  $\beta_j\not=0$ for all $j$, \eqref{cor:picewisfunctionsPolynomialConst} is equivalent to the following: for $t\in G_i\cap H_j$, $1\leq i\leq l$, $1\leq j\leq m$, $(\alpha_i a(t)-\delta_0)b(t)=0$, which is equivalent to that the set
  \begin{align*}
  \supp (\alpha_ia(t)-\delta_0)\cap \supp b(t)\cap (G_i\cap H_j)
  \end{align*}
  has measure zero, for all $i\in \{1,\ldots, l\}$, $j\in\{1,\ldots,m\}$.
   \qed
  \end{proof}

In the following corollary we consider a case of the operator $B$ corresponding to one of the subsets of the partition
defining the operator $A$.
  \begin{cor:DjinjaTumwSilv1}\label{corPiewiseOpSamePartPropertyLikeProjection}
Let $(\mathbb{R},\Sigma,\tilde{m})$ the standard Lebesgue measure space in the real line.
Let $A:L_p(\mathbb{R})\to L_p(\mathbb{R})$,  $B:L_p(\mathbb{R})\to L_p(\mathbb{R})$, $1\le p\le\infty$ be defined as follows, for almost every $t$,
\begin{gather*}
  Ax(t)=\sum\limits_{i=1}^{l} \alpha_i a(t) I_{G_i}(t) x(t ), \qquad (Bx)(t)= I_{G_k}(t)b(t)x(t), \quad 1\le k\le l
\end{gather*}
where $G_i\in \Sigma$ for each $i=1,\ldots,l$, $\tilde{m}( G_i \cap G_j)=0$ for $i\not=j$, $I_{G_k}(t)$ is the indicator function for the set $G_k$, $a,b:\mathbb{R}\to\mathbb{R} $ are measurable functions.
   For a polynomial $F(z)=\sum\limits_{j=1}^{n}\delta_j z^j$ with
$\delta_j\in \mathbb{R}$ for $j=1,\dots,n$, the commutation relation $AB=BF(A)$ is satisfied if and only if
\begin{gather*}
  \supp  b\cap \supp \left(\alpha_k a-F(\alpha_k a)\right)\cap G_k
\end{gather*}
has measure zero.
\end{cor:DjinjaTumwSilv1}

\begin{proof}
We write the operator $B$ as follows
\begin{gather*}
  (Bx)(t)= I_{G_k}(t)b(t)x(t)=\sum_{j=1}^{l} \beta_j b(t)I_{G_j}(t)x(t),
\end{gather*}
almost everywhere, where
$\beta_j=\left\{\begin{array}{l}
                  \, 1,\quad j=k \\
                  \, 0,\  \mbox{ otherwise }
                 \end{array}
                 \right. .$
 By using Corollary \ref{corSamePartitionPiecwiseFunct} we have $AB=BF(A)$ if and only if, for almost very $t$,
\begin{gather*}
  \sum\limits_{i=1}^l \alpha_i  \beta_i a(t)b(t)I_{G_i}(t) =\sum\limits_{i=1}^l b(t) F(a(t)\alpha_i )\beta_i I_{G_i}(t).
\end{gather*}
We can simplify the last relation as follows
\begin{gather*}
  \alpha_k a(t)b(t) I_{G_k}(t)=b(t)F(a(t)\alpha_k) I_{G_k}(t) \
  \Leftrightarrow \ b(t)(a(t) \alpha_k-F(a(t)\alpha_k))=0,
\end{gather*}
for almost every $t\in G_k$, which is equivalent to that the set
\begin{gather*}
  \supp  b\cap \supp  \left(\alpha_k a-F(\alpha_k a)\right)\cap G_k
\end{gather*}
has measure zero.
 \qed
\end{proof}

\begin{ex:DjinjaTumwSilv1}{\rm
Let $A:L_p[0,1]\to L_p[0,1]$,  $B:L_p[0,1]\to L_p[0,1]$, $1\le p\le \infty$ be defined as follows
\begin{align*}
  & (Ax)(t)=\alpha I_{[0,1/3]}(t)x(t)+\beta I_{]1/3,1/2[}(t)x(t)+\gamma I_{]1/2,1]}(t)x(t),\\
  & (Bx)(t)= I_{]1/3, 1/2[}(t)x(t),
\end{align*}
where $\alpha,\, \beta$ are constants and $I_E$ is the indicator function of the set $E$.
For a monomial $F(z)=z^n$, where $n$ is a positive integer, we have
\begin{equation*}
  AB=BF(A)\quad (\mbox{that is}\ AB=BA^n)
\end{equation*}
if and only if
$F(\beta)=\beta.$
In fact, taking a partition $G_1\cup G_2\cup G_3$ where
\begin{equation*}
  G_1=[0,\frac{1}{3}],\quad G_2=]\frac{1}{3},\frac{1}{2}[,\quad G_3=[\frac{1}{2},1],
\end{equation*}
we have  $(Ax)(t)=(\alpha I_{G_1}(t)+\beta I_{G_2}(t) + \gamma I_{G_3}(t))x(t)$ and $(Bx)(t)=I_{G_2}(t)x(t)$. Note that $a(\cdot)=b(\cdot)=1$. By  Corollary \ref{corPiewiseOpSamePartPropertyLikeProjection}, we have $AB=BF(A)$ if and only if
$$ \supp  b\cap \supp  \left(\alpha_k a-F(\alpha_k a)\right)\cap G_k=]\frac{1}{3},\frac{1}{2}[\ \cap\ \supp (\beta-F(\beta))$$
has measure zero which holds if and only if $F(\beta)=\beta$.

In this case, we  can also get the same result by the following direct computation. For almost every $t\in [0,1]$,
\begin{align*}
  &(ABx)(t)=A(Bx)(t)=(\alpha I_{[0,1/3]}(t)+ \beta I_{]1/3,1/2[}(t) +\gamma I_{]1/2,1]}(t) )(Bx)(t) \\
  &=(\alpha I_{[0,1/3]}(t)+\beta I_{]1/3,1/2[}(t)+\gamma I_{]1/2,1]}(t) )I_{]1/3,1/2[(t)}x(t)=\beta I_{]1/3,1/2[(t)}x(t) \\
  &(A^2x)(t)=A(Ax)(t)=(\alpha I_{[0,1/3]}+\beta I_{]1/3,1/2[}+\gamma I_{]1/2,1]} )(Ax)(t)
  \\
  &=(\alpha I_{[0,1/3]}(t)+\beta I_{]1/3,1/2[}(t)+\gamma I_{]1/2,1]}(t) )
  \\
  &
    \hspace{2cm}\cdot (\alpha I_{[0,1/3]}(t)+\beta I_{]1/3,1/2[}(t)+\gamma I_{]1/2,1]}(t) )x(t)\\
  &=(\alpha^2 I_{[0,1/3]}(t)+\beta^2I_{]1/3,1/2[}(t)+\gamma^2 I_{]1/2,1]}(t) )x(t).
  \end{align*}
  In general, for almost every $t$,
  \begin{gather*}
  (A^nx)(t)=(\alpha^n I_{[0,1/3]}(t)+\beta^n I_{]1/3,1/2[}(t) +\gamma^n I_{[1/2,1]}(t))x(t).
  \end{gather*}
  Thus, for almost every $t$,
  \begin{eqnarray*}
  (BF(A)x)(t)&=&(BA^nx)(t)=B(A^nx)(t)= I_{]1/3,1/2[}(t)(A^nx)(t)\\
             &=&\beta^n I_{]1/3,1/2[}(t)x(t)=F(\beta)I_{]1/3,1/2[}(t)x(t).
\end{eqnarray*}
}
\end{ex:DjinjaTumwSilv1}



\section{Representations involving inner superposition operators }\label{SecInnerSuperOpDST2}

In this section we will look at the commutation relation $BA=F(A)B$ which has important applications in wavelets analysis for instance.
When $B=0$, the relation $BA=F(A)B$ is trivially satisfied for any $A$.
If $A=0$ then the relation $BA=F(A)B$ reduces to $F(0)B=0$. This implies either $B=0$, or
$F(0)=0$ with $B$ being any well defined operator. Thus, we focus on construction and properties of non-zero
representations of $BA=F(A)B$.

\begin{prop:DjinjaTumwSilv1}\label{PropInnerSuperpositionOpWavelets}
Let $(\mathbb{R},\Sigma,\tilde{m})$ the standard Lebesgue measure space in the real line.
Let $A:L_p(\mathbb{R})\to L_p(\mathbb{R})$,  $B:L_p(\mathbb{R})\to L_p(\mathbb{R})$, $ p\ge 1$ be defined as follows, for almost every $t\in\mathbb{R}$,
\begin{gather*}
  (Ax)(t)=\sum\limits_{i=-\infty}^{\infty} \alpha_i  I_{G_i}(t) x(t-1 ), \qquad (Bx)(t)= \beta x(\gamma t),
\end{gather*}
where $G_i\in \Sigma$ for each $i\in\mathbb{Z}$, $\tilde{m}(G_i \cap G_j)=0$  for $i\not=j$,  $I_{G_k}(t)$ is the indicator function for the set $G_k$, 
$\beta,\gamma\in\mathbb{R}\setminus\{0\} $  and $\gamma >0$, $\alpha_i\in \mathbb{R}$, for all $i\in \mathbb{Z}$.
For a monomial $F(z)=\delta z^m$, where $\delta \in\mathbb{R}\setminus\{0\}$ and $m>0$ is a positive integer,
$(BAx)(t)=(F(A)Bx)(t)$ is satisfied for almost every $t$ if and only if
\begin{gather}\label{PropWaveletsExCondOnePeriodicity}
 \sum\limits_{i=- \infty}^{\infty}\alpha_i I_{G_i}(\gamma t)x(\gamma t-1)=\delta\sum\limits_{i=- \infty}^{\infty}\alpha^m_i I_{G_i}(t)x(\gamma t-\gamma m),
\end{gather}
for almost every $t$.
\end{prop:DjinjaTumwSilv1}

\begin{proof}
 For almost every $t\in \mathbb{R}$ we have
\begin{eqnarray*}
  (A^2x)(t)&=&A\left(\sum\limits_{i=-\infty}^{\infty} \alpha_i  I_{G_i}(t) x(t-1 )\right)=\sum\limits_{i=-\infty}^{\infty} \alpha_i  I_{G_i}(t) (Ax)(t-1 )\\
  &=&\sum\limits_{i=-\infty}^{\infty} \alpha_i  I_{G_i}(t) \sum\limits_{j=-\infty}^{\infty} \alpha_j  I_{G_j}(t) x(t-2 )=
  \sum\limits_{i=-\infty}^{\infty} \alpha^2_i  I_{G_i}(t) x(t-2 ).
  \end{eqnarray*}
  We suppose that for almost every $t\in\mathbb{R}$,
  \begin{gather*}
  (A^mx)(t)=\sum\limits_{i=-\infty}^{\infty} \alpha^m_i  I_{G_i}(t) x(t-m ),\quad  \text{for}\ m\geq 1
  \end{gather*}
 Then we have, for almost every $t$,
 \begin{eqnarray*}
  (A^{m+1}x)(t)&=&A(A^mx)(t)=A\left(\sum\limits_{i=-\infty}^{\infty} \alpha^m_i  I_{G_i}(t) x(t-m )\right)\\
  &=& \sum\limits_{i=-\infty}^{\infty} \alpha_i   I_{G_i}(t) \sum\limits_{j=-\infty}^{\infty} \alpha^m_j  I_{G_j}(t) x(t-m-1 )\\
  &=&\sum\limits_{i=-\infty}^{\infty} \alpha^{m+1}_i  I_{G_i}(t) x(t-(m+1)).
  \end{eqnarray*}
Moreover, we have for almost every $t\in\mathbb{R}$,
  \begin{eqnarray*}
  (BAx)(t)&=&
  \beta \sum\limits_{i=-\infty}^{\infty} \alpha_i  I_{ G_i}(\gamma t) x(\gamma t-1) \\
  \delta(A^m Bx)(t)&=&
  \delta\beta \sum\limits_{i=- \infty}^{\infty} \alpha^m_i  I_{G_i}(t) x(\gamma t- \gamma m ).
\end{eqnarray*}
Therefore, $(BAx)(t)=(F(A)Bx)(t)$ for almost every $t$ if and only if
\begin{align}
& \sum\limits_{i=- \infty}^{\infty}\alpha_i I_{G_i}(\gamma t)x(\gamma t-1)=\delta\sum\limits_{i=- \infty}^{\infty}\alpha^m_i I_{G_i}(t)x(\gamma t-\gamma m), \nonumber \\
&\hspace{4cm} \mbox{for almost every}\ t, \ i\in\mathbb{Z}.  \tag*{\qed}
\end{align}
\end{proof}
\begin{prop:DjinjaTumwSilv1}\label{PropInnerSuperpositionOperatorWavelets}
Let $A:L_p(\mathbb{R})\to L_p(\mathbb{R})$,  $B:L_p(\mathbb{R})\to L_p(\mathbb{R})$, $1\le p \le \infty$ be defined as follows,
$(Ax)(\cdot)=\alpha x(\upsilon(\cdot)),$ $(Bx)(\cdot)= \beta x(\sigma(\cdot)),
$
where $\upsilon,\sigma:\mathbb{R}\to\mathbb{R}$ are continuous functions and $\alpha, \beta \in \mathbb{R}\setminus \{0\}$.
For a  monomial   $F(z)=\delta z^m$ where $\delta\in\mathbb{R}\setminus\{0\}$,
$m\in \mathbb{Z}_{>0}=\{1,2,\ldots\}$, for $x\in L_p(\mathbb{R})$ and for almost every $t$,
\begin{gather*}
  (BAx)(t)=(F(A)Bx)(t),
\end{gather*}
if and only if
$ x(\upsilon(\sigma(t)))=\delta\alpha^{m-1}  x(\sigma(\upsilon^{\circ(m)})(t))$.
\end{prop:DjinjaTumwSilv1}

\begin{proof}
  We have for almost every $t$
  \begin{gather*}
    (BAx)(t)=\alpha \beta x(\upsilon(\sigma(t))), \\
    (A^2x)(t)=\alpha^2 x(\upsilon(\upsilon(t))).
    \end{gather*}
  In the same way we have for almost every $t$
    \begin{gather*}
    (A^mx)(t)=\alpha^m x(\upsilon^{\circ(m)}(t)),\\
    \delta(A^mBx)(t)=\delta\alpha^m \beta x(\sigma(\upsilon^{\circ(m)})(t)).
  \end{gather*}
  Then, for all $L_p(\mathbb{R})$, $1\leq p \leq \infty$, $BAx=(F(A)B)x$  if and only if,  for almost every $t$,
  \begin{equation} 
    \alpha\beta x(\upsilon(\sigma(t)))=\delta\alpha^m \beta x(\sigma(\upsilon^{\circ(m)})(t)). \tag*{\qed}
  \end{equation}
\end{proof}

\begin{ex:DjinjaTumwSilv1}{\rm
Consider  operators $A:L_2(\mathbb{R})\to L_2(\mathbb{R})$, $B:L_2(\mathbb{R})\to L_2(\mathbb{R})$  be defined as follows, for almost every $t$,
$$
(Ax)(t)=x(t-1),\quad  (Bx)(t)=\gamma^{1/2}x(\gamma t),\quad \gamma>0.
$$
These operators are
particular case of the corresponding ones of Proposition \ref{PropInnerSuperpositionOperatorWavelets} when
$\upsilon(t)=t-1$, $\sigma(t)={\gamma t}$, $\alpha=1$, $\beta=\gamma^{\frac{1}{2}}$.
We get
$$
\upsilon(\sigma(t))={\gamma t-1},\quad \sigma(\upsilon^{\circ(m)})(t)=\gamma t-\gamma m,\quad m\in \mathbb{Z}_{>0}=\{1,2,\ldots\}$$
Then $BAx=\delta A^mBx$ for $x\in L_p(\mathbb{R})$, for almost every $t$, if and only
$$ x(\gamma t-1)=\delta x(\gamma t-\gamma m).$$
If $\gamma=\frac{1}{m}$ and $\delta=1$, then we get on both sides $x(\frac{t}{m}-1)$ and thus  $BA=A^m B$ on $L_p(\mathbb{R})$.
For example, when $\gamma=1/2$ and $m=2$, we have operators well known operators important for example in wavelets theory,
\begin{equation*}
  (Ax)(t)=x\left(t-1\right),\quad  (Bx)(t)=\frac{1}{2^{1/2}}x\left(\frac{1}{2}t\right),
\end{equation*}
 satisfying $  BA=A^2B$.
}
\end{ex:DjinjaTumwSilv1}

\section{Representations involving weighted composition operators}\label{SecMultiplOpDST2}
In this section we consider pairs of operators $(A,B)$ which involve weighted composition operators, and the multiplication operator composed with the point evaluation functional $\upsilon_{\gamma}:x(t)\mapsto x(\gamma)$  ("boundary value" operator, "one-dimensional range" operator).  These operators act on some spaces of real-valued functions $x(t)$ of one real variable by the formulas
$$
(T_{w,\sigma} x)(t)= w(t) x(\sigma(t)), \quad T_{w,\upsilon_{\gamma}} x(t)= w(t) x(\gamma), \quad t, \gamma \in \mathbb{R}.
$$
Note that the multiplication operator composed with the point evaluation functional is a special case of weighted composition operator in these notations since
$$T_{w,\upsilon_{\gamma}}=T_{w,\sigma},\quad \mbox{for}\quad  \sigma=\upsilon_{\gamma}.$$

\begin{prop:DjinjaTumwSilv1}\label{propDiagOperatorLp}
Let $A:L_p(\mathbb{R})\to L_p(\mathbb{R})$,  $B:L_p(\mathbb{R})\to L_p(\mathbb{R})$, $1\le p \le \infty$ be defined for measurable functions $a,b:\mathbb{R}\to\mathbb{R} $ by
\begin{equation*}
  (Ax)(\cdot)=a(\cdot) x(\cdot),\quad (Bx)(\cdot)= b(\cdot) x(\cdot).
\end{equation*}
For a polynomial $F(z)=\sum\limits_{j=0}^{n}\delta_j z^j$
with $\delta_0,\ldots,\delta_n\in \mathbb{R}$,
\begin{equation*}
  AB=BF(A),
\end{equation*}
if and only if the set $\supp b\ \cap \ \supp (a-F(a))$
has measure zero.
\end{prop:DjinjaTumwSilv1}

\begin{proof}
For almost every $t\in \mathbb{R}$,
\begin{gather*}
  (ABx)(t)=A(Bx)(t)=a(t) b(t) x(t), \\
  (A^2x)(t)=A(Ax)(t)=a(t) (Ax)(t) =a(t)^2 x(t)=a(t)^2 x(t),\\
  (A^3x)(t)=A(A^2x)(t)=a(t)^3 x(t).
  \end{gather*}
 Therefore, for $m\ge 1$ and for almost every $t$,
 \begin{gather*}
  (A^mx)(t)=a(t)^m x(t),\
  (BA^mx)(t)=B(A^mx)(t)= b(t) (A^mx)(t)=b(t)a(t)^m x(t).
  \end{gather*}

Thus, we have for almost every $t$,
\begin{equation*}
 (BF(A)x)(t)=\sum_{k=0}^{n} \delta_k(B(A^k x))(t)=\sum_{k=0}^{n} \delta_k b(t)a(t)^k x(t)=b(t)F(a(t))x(t).
\end{equation*}
Then  $AB=BF(A)$ if and only if for almost every $t$,
\begin{equation*}
  b(t)F(a(t))=a(t)b(t).
\end{equation*}
 This is equivalent to the set
$
  \supp b\ \cap \ \supp (a-F(a))
$
having measure zero. \qed
\end{proof}

\subsection{Representations by operators on \texorpdfstring{$C[\alpha,\beta]$}{C[\alpha,\beta]}}
In this subsection we consider pairs of operators $(A,B)$ which involve weighted composition operators, and the multiplication operator composed with the point evaluation functional $\upsilon_{\gamma}:x(t)\mapsto x(\gamma)$
("boundary value" operator, "one-dimensional range" operator).  These operators act on some spaces of continuous real-valued functions of one real variable $x(t)$ by the formulas:
$$
(T_{w,\sigma} x)(t)= w(t) x(\sigma(t)), \quad T_{w,\upsilon_{\gamma}} x(t)= w(t) x(\gamma), \quad t, \gamma \in \mathbb{R}.
$$

\begin{prop:DjinjaTumwSilv1}
Let $A:C[\alpha,\beta]\to C[\alpha,\beta]$,  $B:C[\alpha,\beta]\to C[\alpha,\beta]$, $1\le p \le \infty$ defined as follows
\begin{equation} \label{ops:weightedcomposops}
  (Ax)(t)=a(t) x(\upsilon(t)),\quad (Bx)(t)= b(t) x(\sigma(t)),
\end{equation}
where $\alpha, \beta$ are real numbers, $\alpha<\beta$,  $a,b,\upsilon, \sigma :[\alpha,\beta]\to[\alpha,\beta]$ are continuous functions.
For a  monomial   $F(z)=\delta z^m$ where $\delta\in\mathbb{R}\setminus\{0\}$, $m\in \mathbb{Z}_{>0}=\{1,2,\ldots\}$, and for $x\in C[\alpha,\beta]$ and for each $t\in [\alpha,\beta]$,
\begin{equation*}
  (ABx)(t)=(BF(A)x)(t),
\end{equation*}
if and only if
\begin{multline*}
   a(t)b(\upsilon(t))x(\sigma(\upsilon(t)))\\
   =\delta b(t) a(\sigma(t))\cdot a(\upsilon(\sigma (t)))\cdot\ldots\cdot a(\upsilon^{\circ (m-1)}(\sigma (t))) x(\upsilon^{\circ(m)}(\sigma(t))).
  \end{multline*}
In particular, if $x\in C[\alpha,\beta]$ satisfies
  $x(\sigma(\upsilon (\cdot))))=x(\upsilon^{\circ(m)}(\sigma(\cdot))))$, then
  $ABx=BF(A)x$ if and only if
  \begin{equation*}
   a(t)b(\upsilon(t))=\delta b(t) a(\sigma(t))\cdot a(\upsilon(\sigma (t)))\cdot\ldots\cdot a(\upsilon^{\circ (m-1)}(\sigma (t)))
  \end{equation*}
for all $t\in \supp x(\sigma(\upsilon (\cdot))))\bigcap [\alpha,\beta]$
  \end{prop:DjinjaTumwSilv1}

\begin{proof}
  For each $t \in [\alpha,\beta]$,
  \begin{eqnarray*}
    (ABx)(t)&=&a(t)b(\upsilon(t))x(\sigma(\upsilon(t))), \\
    (A^2x)(t)&=&a(t)a(\upsilon(t)) x(\upsilon(\upsilon(t)))).
    \end{eqnarray*}
  In the same way we have for $n\ge 1$ and for each $t\in [\alpha,\beta]$,
    \begin{gather*}
    (A^nx)(t)=a(t)a(\upsilon(t))a(\upsilon^{\circ(2)}(t))\ldots a(\upsilon^{\circ(n-1)}(t)) x(\upsilon^{\circ(n)}(t)),\\
    \delta(BA^m x)(t)=\delta b(t) a(\sigma(t))a(\upsilon(\sigma(t)))\cdot\ldots \cdot a(\upsilon^{\circ(m-1)}(\sigma(t))) x(\upsilon^{\circ(m)}(\sigma(t))).
  \end{gather*}
  Then $ABx=BF(A)x$ for all $x\in C[\alpha,\beta]$ if and only if for all $t\in [\alpha,\beta]$,
  \begin{multline*}
   a(t)b(\upsilon(t))x(\sigma(\upsilon(t))) \\ =\delta b(t) a(\sigma(t))\cdot a(\upsilon(\sigma(t)))\cdot\ldots\cdot a(\upsilon^{\circ(m-1)}(\sigma(t)))  x(\upsilon^{\circ(m)}(\sigma(t))).
  \end{multline*}
  If $x(\sigma(\upsilon(\cdot)))=x(\upsilon^{\circ(m)}(\sigma(\cdot)))$, then $ABx=BF(A)x$
  if and only if
  \begin{gather*}
  a(t)b(\upsilon(t))=\delta b(t) a(\sigma(t))\cdot a(\upsilon(\sigma(t)))\cdot\ldots\cdot a(\upsilon^{\circ(m-1)}(\sigma(t))).
  \end{gather*}
for $t\in [\alpha,\beta]$ such that $x(\sigma(\upsilon(t)))\neq 0$, since $(ABx)(t)=0=(BF(A)x)(t)$ for $t\in [\alpha,\beta]$
such that $x(\sigma(\upsilon(t)))=0$. \qed
\end{proof}

\subsubsection{Representations when \texorpdfstring{$B$}{B} is the multiplication operator}
Let us consider the operators \eqref{ops:weightedcomposops}
with constant $\upsilon(\cdot)=\gamma$ and with $\sigma_B(t)=t, t\in \mathbb{R}$.

\begin{lem:DjinjaTumwSilv1}\label{LemmaContinuousFunctionsSupportMeasureZero}
Let $a,b\in C[\alpha,\beta]$, $\alpha,\beta\in\mathbb{R}$, $\alpha<\beta$.
If the set $\Omega=\supp a\, \cap \, \supp b$ has Lebesgue measure zero, then it is empty.
\end{lem:DjinjaTumwSilv1}

\begin{proof}
Without loss of generality we suppose that there exists $\alpha_0\in ]\alpha,\beta[$ such that
$\alpha_0\in \Omega$, that is, $\Omega$ is not empty. Then, $a(\alpha_0)\not=0$ and $b(\alpha_0)\not=0$.
Since $a,b$ are continuous functions, there is an open interval $V_{\alpha_0}=]\alpha_0-\varepsilon,\alpha_0+\varepsilon[\ \subseteq\ ]\alpha,\beta[$, $\varepsilon>0$ such that $a(t)\not=0$ and $b(t)\not=0$
for all $t\in V_{\alpha_0}$.
Then  the Lebesgue measure of $\Omega$ is positive. But this contradicts the hypothesis. Then the set $\Omega$ must be empty. \qed
\end{proof}

\begin{prop:DjinjaTumwSilv1} \label{PropContFunctAtPoint}
Let $A:C[\alpha,\beta]\to C[\alpha,\beta]$,  $B:C[\alpha,\beta]\to C[\alpha,\beta]$ defined by
\begin{equation*}
  (Ax)(t)=a(t) x(\gamma),\quad (Bx)(t)= b(t) x(t),
\end{equation*}
where $\alpha, \beta$ are real numbers, $\alpha<\beta$, $\gamma\in [\alpha,\beta]$   and $a,b:[\alpha,\beta]\to\mathbb{R}$ are continuous functions.
Let $F(z)=\sum\limits_{j=0}^{n}\delta_j z^j$,  $\delta_j\in\mathbb{R} $, $j=0,\ldots,n$.
For $x\in C[\alpha,\beta]$ and $t\in[\alpha,\beta] $,
\begin{equation*}
  (ABx)(t)=(BF(A)x)(t)
\end{equation*}
if and only if,
\begin{equation}
   a(t)b(\gamma)x(\gamma)=\delta_0  b(t)x(t)+
   b(t) a(t)x(\gamma)\sum\limits_{j=1}^{n} \delta_j  a(\gamma)^{j-1}.
\end{equation}
\end{prop:DjinjaTumwSilv1}

\begin{proof}
We have
\begin{gather*}
  (ABx)(t)=A(Bx)(t)=a(t)(Bx)(\gamma)=a(t)b(\gamma)x(\gamma) \\
  (A^mx)(t)= a(t) a(\gamma)^{m-1} x(\gamma),\quad m=1,2,\ldots\\
  ((BA^m)x)(t)=(BA^mx)(t)=
   b(t) a(t) a(\gamma)^{m-1} x(\gamma), \ m=1,2,\ldots\\
  (BF(A)x)(t)= \delta_0  b(t)x(t)+\sum\limits_{j=1}^{n} \delta_j b(t) a(t) a(\gamma)^{j-1} x(\gamma).
\end{gather*}
Then $(ABx)(t)=(BF(A)x)(t)$ if and only if
\begin{equation*}
  a(t)b(\gamma)x(\gamma)=\delta_0  b(t)x(t)+\sum\limits_{j=1}^{n} \delta_j b(t) a(t) a(\gamma)^{j-1} x(\gamma). \ \qed
\end{equation*}
\end{proof}

\begin{cor:DjinjaTumwSilv1} \label{CorContFunctAtPointFzerozero}
Let $A:C[\alpha,\beta]\to C[\alpha,\beta]$,  $B:C[\alpha,\beta]\to C[\alpha,\beta]$ be defined by
\begin{equation*}
  (Ax)(t)=a(t) x(\gamma),\quad (Bx)(t)= b(t) x(t),
\end{equation*}
where $\alpha, \beta$ are real numbers, $\alpha<\beta$, $\gamma\in [\alpha,\beta]$   and $a,b:[\alpha,\beta]\to\mathbb{R}$ are continuous functions
such that $a(\gamma)\not=0$, $b(\gamma)\not=0$.
Let $F(z)=\sum\limits_{j=1}^{n}\delta_j z^j$,  where $\delta_j\in \mathbb{R}$, $j=1,\ldots,n$.
The following statements hold.
\begin{enumerate}
\item If $b(\cdot)$ is not  constant in any open interval included in $[\alpha,\beta]$, then
\begin{equation*}
  AB\not=BF(A).
\end{equation*}
\item If $b(\cdot)$ is  constant in some open interval included in $[\alpha,\beta]$, then
\begin{equation*}
  AB=BF(A)
\end{equation*}
if and only if
$
  \supp a\,\cap \, \supp\left(b(\gamma)- b\cdot k_1\right)=\emptyset,
$
 where
$  k_1=\sum\limits_{j=1}^{n} \delta_j a(\gamma)^{j-1}.$

In particular, if $b(\cdot)=b(\gamma)$ is identically non-zero constant in $[\alpha,\beta]$ then $AB=BF(A)$ if and only if $k_1=1$.
\end{enumerate}
\end{cor:DjinjaTumwSilv1}

\begin{proof}
Suppose that there exist two continuous functions $a,\ b$, $a(\gamma)\not=0$, $b(\gamma)\not=0$ and a polynomial
$F(z)=\sum\limits_{j=1}^{n}\delta_j z^j$ such that $ABx=BF(A)x$ for all $x\in C[\alpha,\beta]$.
 By applying Proposition \ref{PropContFunctAtPoint} we have $AB=BF(A)$  if and only  if for all $x\in C[\alpha,\beta]$
 \begin{equation*}
     a(t)b(\gamma)x(\gamma)=\sum\limits_{j=1}^{n} \delta_j b(t) a(t) a(\gamma)^{j-1} x(\gamma).
 \end{equation*}
 This is equivalent to
 \begin{equation*}
   a(t)b(\gamma)=\sum\limits_{j=1}^{n} \delta_j b(t) a(t) a(\gamma)^{j-1} \Leftrightarrow
   a(t)\left(b(\gamma)- b(t)\sum\limits_{j=1}^{n} \delta_j a(\gamma)^{j-1} \right)=0,
\end{equation*}
 which is equivalent to the set
  \begin{equation*}
\Omega= \supp a\,\cap \, \supp\left(b(\gamma)- b\cdot\sum\limits_{j=1}^{n} \delta_j a(\gamma)^{j-1} \right)
\end{equation*}
has measure zero. Since $a,b$ are  continuous functions, by applying Lemma \ref{LemmaContinuousFunctionsSupportMeasureZero} we conclude that the set $\Omega$ is empty.
We consider the following cases.
\begin{enumerate}[label=,leftmargin=*]
\item {\bf Case 1:} If $b(\cdot)$ is not constant in any open interval included in $[\alpha,\beta]$ and
\begin{equation*}
  k_1=\sum\limits_{j=1}^{n} \delta_j a(\gamma)^{j-1}\not=1,
\end{equation*}
then $\gamma\in \Omega$. This contradicts the fact that $\Omega$ is empty. \par
\item {\bf Case 2:} Suppose $k_1=1$ and $b(\cdot)$ is not  constant in any open interval included
in $[\alpha,\beta]$ and without loss of generality suppose that $\gamma \in ]\alpha,\beta[$.
Since the functions $a$ and $b(\gamma)-k_1b(t)=b(\gamma)-b(t)$ are continuous,  then for some positive $\varepsilon$,
 we can find an open interval
$]\gamma-\varepsilon, \ \gamma+\varepsilon [\subseteq ]\alpha,\beta[$  such that the set
\begin{equation*}
  \supp a\,\cap \, \supp\left(b(\gamma)- b\right)\, \cap\, ]\gamma-\varepsilon, \ \gamma+\varepsilon [ \not=\emptyset,
\end{equation*}
which is a contradiction. \par
\item {\bf Case 3:} If $b(\cdot)$ is constant in an open interval contained in $[\alpha,\beta]$ then $AB=BF(A)$ if and only if
  \begin{equation*}
\Omega= \supp a\,\cap \, \supp\left(b(\gamma)- b\cdot\sum\limits_{j=1}^{n} \delta_j a(\gamma)^{j-1} \right)=\emptyset.
\end{equation*}
\item {\bf Case 4: } If  $b(\cdot)$ is identically constant in $[\alpha,\beta]$ and $k_1=1$ then $\Omega$ is empty. If $k_1\not=1$ then $\Omega=\supp  a $ must be empty. This implies $a(\cdot)\equiv0$, but, $a(\gamma)\not= 0$. That is a contradiction.
    Therefore, the condition of Proposition
\ref{PropContFunctAtPoint} is fulfilled. \qed
\end{enumerate}
\end{proof}

\begin{rem:DjinjaTumwSilv1}
Note that if operators  $A$ and $B$ defined in Corollary \ref{CorContFunctAtPointFzerozero} are such that $AB=BF(A)$ for some polynomial
$F(z)=\sum\limits_{j=1}^{n}\delta_j z^j,$ then $k_1=\sum\limits_{j=1}^{n}\delta_j a(\gamma)^{j-1}=1$. In fact, if $k_1\not=1$, then $\gamma \in \Omega$, thus by Proposition \ref{PropContFunctAtPoint}, $AB\not=BF(A)$.

\end{rem:DjinjaTumwSilv1}

\begin{ex:DjinjaTumwSilv1}{\rm
Let $A:C[0,3]\to C[0,3]$,  $B:C[0,3]\to C[0,3]$ be defined as follows
\begin{equation*}
  (Ax)(t)=a(t) x(1),\quad (Bx)(t)= b(t) x(t),
\end{equation*}
where
$a(t) = \left\{ \begin{array}{rr} t(2-t), & \mbox{ if }\ 0\le t\le 2 \\
     0, & \mbox{ otherwise, } \end{array}\right.$ and
$ b(t) = \left\{ \begin{array}{rr} 1, & \mbox{ if }\ 0\le t\le 2 \\
     (t-1)^3, & \mbox{ if }\ 2<t\le 3 \end{array}.\right.$
Consider a polynomial
 $F(z)=\sum\limits_{j=1}^{n}\delta_j z^j$, $\delta_j\in\mathbb{R}$, $j=1,\ldots,n$. If
 $
 k_1=\sum\limits_{j=1}^{n} \delta_j=1, 
 $
 then these operators satisfy $AB=BF(A)$. Indeed, by  Corollary \ref{CorContFunctAtPointFzerozero}
 we have
 \begin{equation*}
   \Omega=\supp a\, \cap \, (b(1)-b)=\emptyset.
 \end{equation*}
}\end{ex:DjinjaTumwSilv1}

\begin{cor:DjinjaTumwSilv1}\label{corRepreContFunctBmultOpagammazero}
Let $A:C[\alpha,\beta]\to C[\alpha,\beta]$,  $B:C[\alpha,\beta]\to C[\alpha,\beta]$ be defined by
\begin{equation*}
  (Ax)(t)=a(t) x(\gamma),\quad (Bx)(t)= b(t) x(t),
\end{equation*}
where $\alpha, \beta$ are real numbers, $\alpha<\beta$, $\gamma\in [\alpha,\beta]$   and $a,b:[\alpha,\beta]\to\mathbb{R}$ are continuous functions such that $b(\gamma)\not=0$.
Consider a polynomial $F(z)=\sum\limits_{j=0}^{n}\delta_j z^j$ with $\delta_j\in\mathbb{R}$ for $j=0,\ldots,n$. If $a(\gamma)=0$, then
$  AB=BF(A)$ if and only if $\delta_0=0$ and
\begin{equation*}
      \Omega= \supp a \, \cap \, \supp (b(\gamma)-\delta_1 b) = \emptyset.
    \end{equation*}
\end{cor:DjinjaTumwSilv1}

\begin{proof}
 By applying Proposition \ref{PropContFunctAtPoint} we have $(ABx)(t)=(BF(A)x)(t)$  if and only  if, for all $x\in C[\alpha,\beta]$,
 \begin{equation*}
     a(t)b(\gamma)x(\gamma)=\delta_0 b(t)x(t)+\delta_1 b(t) a(t)x(\gamma) +\sum\limits_{j=2}^{n} \delta_j b(t) a(t) a(\gamma)^{j-1} x(\gamma).
 \end{equation*}
  By taking $t=\gamma$ and if $a(\gamma)=0$ we get $\delta_0 b(\gamma)x(\gamma)=0$ for all $x\in C[\alpha,\beta]$. Since
  $b(\gamma)\not=0$ then $\delta_0$ must be zero. By using this, we remain with the equation
  \begin{equation*}
     a(\cdot)b(\gamma)x(\gamma)=\delta_1 a(\cdot)b(\cdot)x(\gamma)
  \end{equation*}
  for all $x\in C[\alpha,\beta]$.  This is equivalent to
    $ a(\cdot)(b(\gamma)-\delta_1 b(\cdot))=0,
    $
    which itself is equivalent to the set
    $
      \Omega= \supp a \, \cap \, \supp (b(\gamma)-\delta_1 b)
    $
    having measure zero. Since the functions involved are continuous then the set  $\Omega$ is empty.
    \qed
\end{proof}

\begin{rem:DjinjaTumwSilv1}{\rm
According to Corollary \ref{corRepreContFunctBmultOpagammazero}, if $a(\gamma)=0$  the role of the polynomial $F(\cdot)$ is only played by the coefficient $\delta_1$. In particular, Corollary \ref{corRepreContFunctBmultOpagammazero} establishes conditions for representations of the quantum plane relation
$AB=\delta_1 BA$, for a real constant $\delta_1$.}
\end{rem:DjinjaTumwSilv1}

\begin{cor:DjinjaTumwSilv1}\label{corContFunctionsIVOPaGammaNonZerobGammaZero}
Let $A:C[\alpha,\beta]\to C[\alpha,\beta]$,  $B:C[\alpha,\beta]\to C[\alpha,\beta]$ be defined by
\begin{equation*}
  (Ax)(t)=a(t) x(\gamma),\quad (Bx)(t)= b(t) x(t),
\end{equation*}
where $\alpha, \beta$ are real numbers, $\alpha<\beta$, $\gamma\in [\alpha,\beta]$   and $a,b:[\alpha,\beta]\to\mathbb{R}$ are continuous functions  such that $a(\gamma)\not=0$.
Consider a polynomial $F(z)=\sum\limits_{j=1}^{n}\delta_j z^j$ with  $\delta_j\in\mathbb{R}$ for $j=1,\ldots,n$.
Let
$k_1=\sum\limits_{j=1}^{n} \delta_j (a(\gamma))^{j-1}.
$
If $b(\gamma)=0$, then
$  AB=BF(A)$ if and only if
$k_1=0$ or $\supp a \,\cap  \supp b=\emptyset.$
\end{cor:DjinjaTumwSilv1}

\begin{proof}
 By applying Proposition \ref{PropContFunctAtPoint} we have $(ABx)(t)=(BF(A)x)(t)$
 if and only  if for all $x\in C[\alpha,\beta]$ and for all $t\in [\alpha,\beta]$,
 \begin{equation*}
     a(t)b(\gamma)x(\gamma)=\sum\limits_{j=1}^{n} \delta_j b(t) a(t) a(\gamma)^{j-1} x(\gamma).
 \end{equation*}
  By using the hypothesis, this reduces to the following:
  \begin{equation*}
    \forall\, x\in C[\alpha,\beta]:\  k_1 a(t)b(t)x(\gamma)=0.
  \end{equation*}
  This is equivalent to $k_1a(\cdot)b(\cdot)=0$ which is equivalent to $k_1=0$ or to the set
  \begin{equation*}
    \Omega_1=\supp a \, \cap \, \supp b
  \end{equation*}
  having measure zero. Since $a,b$ are continuous functions, the set $\Omega_1$ is empty.
  \qed
  \end{proof}

  \begin{ex:DjinjaTumwSilv1}{\rm
Let $A:C[0,2]\to C[0,2]$,  $B:C[0,2]\to C[0,2]$ be defined as follows
\begin{equation*}
  (Ax)(t)=a(t) x(\frac{1}{2}),\quad (Bx)(t)= b(t) x(t),
\end{equation*}
where
$a(t) = \left\{ \begin{array}{rr} t(1-t), & \mbox{ if }\ 0\le t\le 1 \\
0, & \mbox{ otherwise, } \end{array}\right.$ and
$ b(t) = \left\{ \begin{array}{cc} 0, & \mbox{ if }\ 0\le t\le 1 \\
(1-t)(t-2), & \mbox{ if } 1<t\le 2 \end{array}\right. .
$
These operators satisfy $AB=BF(A)$ for any polynomial
$F(z)=\sum\limits_{j=1}^{n}\delta_j z^j$, $\delta_j\in\mathbb{R}$, $j=1,\ldots,n$. Indeed this follows
 by Corollary \ref{corContFunctionsIVOPaGammaNonZerobGammaZero}.
}\end{ex:DjinjaTumwSilv1}

\begin{ex:DjinjaTumwSilv1}{\rm
Let $A:C[-1,1]\to C[-1,1]$,  $B:C[-1,1]\to C[-1,1]$ be defined by
\begin{equation*}
  (Ax)(t)=a(t) x(0),\quad (Bx)(t)= b(t) x(t),
\end{equation*}
where $a(t)=1+t^2$, $b(t)=t$. Consider a polynomial   $F(z)=\sum\limits_{j=1}^{n}\delta_j z^j$ with
$\delta_j\in\mathbb{R}$, $j=1,\ldots,n$. If $k_1=\delta_1+\ldots+\delta_n=0$,
then these operators satisfy $AB=BF(A)$. Indeed this follows
 by Corollary \ref{corContFunctionsIVOPaGammaNonZerobGammaZero}.
}\end{ex:DjinjaTumwSilv1}

\begin{cor:DjinjaTumwSilv1}\label{corContFunctionsIVOPaGammaZerobGammaZero}
Let $A:C[\alpha,\beta]\to C[\alpha,\beta]$,  $B:C[\alpha,\beta]\to C[\alpha,\beta]$ be defined by
\begin{equation*}
  (Ax)(t)=a(t) x(\gamma),\quad (Bx)(t)= b(t) x(t),
\end{equation*}
where $\alpha, \beta$ are real numbers, $\alpha<\beta$, $\gamma\in [\alpha,\beta]$   and $a,b:[\alpha,\beta]\to\mathbb{R}$ are non-zero continuous functions.
Consider a polynomial $F(z)=\sum\limits_{j=0}^{n}\delta_j z^j$,  $\delta_j\in\mathbb{R}$, $j=0,\ldots,n$. If $a(\gamma)=b(\gamma)=0$ then
$  AB=BF(A)$ if and only if
  $\delta_0=0$ and either $\delta_1=0$ or the set
   $
     \Omega=\supp a \, \cap \, \supp b=\emptyset.
   $
\end{cor:DjinjaTumwSilv1}

\begin{proof}
 By applying Proposition \ref{PropContFunctAtPoint} we have $(ABx)(t)=(BF(A)x)(t)$  if and only  if for all $x\in C[\alpha,\beta]$ and for all $t\in [\alpha,\beta]$,
 \begin{equation*}
     a(t)b(\gamma)x(\gamma)=\delta_0 b(t)x(t)+\delta_1 b(t) a(t) x(\gamma)+\sum\limits_{j=2}^{n} \delta_j b(t) a(t) a(\gamma)^{j-1} x(\gamma).
 \end{equation*}
 By hypothesis, this reduces to the condition
 \begin{equation}\label{CondABBFAproofCorBothagammabgammmazero}
  \forall\, x\in C[\alpha,\beta]:\ \delta_0 b(\cdot)x(\cdot)+\delta_1 a(\cdot)b(\cdot)x(\gamma)=0.
 \end{equation}
 If $\delta_0\not=0$ and $t_0\in \supp b$, then in an open interval
 $]t_0-\varepsilon,t_0+\varepsilon[\subset [\alpha,\beta]$ we have
 $x(t)=-\frac{\delta_1}{\delta_0 }a(t) \zeta, \ t\in ]t_0-\varepsilon,t_0+\varepsilon[,$
  where $\zeta$ is a constant. Since $x(\cdot)$ is continuous in $]t_0-\varepsilon,t_0+\varepsilon[$, thus  $1+x(\cdot)^2$ is also continuous in $]t_0-\varepsilon,t_0+\varepsilon[$, but the identity \eqref{CondABBFAproofCorBothagammabgammmazero} is not valid for this function. This contradicts the condition \eqref{CondABBFAproofCorBothagammabgammmazero}.
  Since $b(\cdot)$ is not identically zero, this implies that $\delta_0=0$. By using this, the condition \eqref{CondABBFAproofCorBothagammabgammmazero}
   reduces to the following condition:
   $
     \forall\, x\in C[\alpha,\beta]:\ \delta_1 a(t)b(t)x(\gamma)=0,
   $
   which is equivalent to $\delta_1=0$ or to
   $
     \Omega=\supp a \, \cap \, \supp b=\emptyset.
   $
   \qed
\end{proof}

\begin{ex:DjinjaTumwSilv1}{\rm
Let $A:C[0,2]\to C[0,2]$,  $B:C[0,2]\to C[0,2]$ be defined by
\begin{equation*}
  (Ax)(t)=a(t) x(1),\quad (Bx)(t)= b(t) x(t),
\end{equation*}
where $a(t)=\sin (\pi t)$, $b(t)=t^2-1$ and $\gamma=1$.
By Corollary \ref{corContFunctionsIVOPaGammaZerobGammaZero}, these operators satisfy $AB=BF(A)$ for any polynomial
 $F(z)=\sum\limits_{j=2}^{n} \delta_j z^j$ with $\delta_j\in\mathbb{R}$, $j=2,\ldots,n$.
}\end{ex:DjinjaTumwSilv1}

\subsubsection{Representations when \texorpdfstring{$A$}{A} is multiplication operator}

\begin{prop:DjinjaTumwSilv1}\label{PropAmultOpBinitialValueOpContFunct}
Let $A:C[\alpha,\beta]\to C[\alpha,\beta]$,  $B:C[\alpha,\beta]\to C[\alpha,\beta]$  be defined by
\begin{equation*}
  (Ax)(t)=a(t)x(t) ,\quad (Bx)(t)= b(t) x(\gamma),
\end{equation*}
where $\alpha,\beta\in\mathbb{R}$, $\alpha<\beta$, $\gamma\in [\alpha,\beta]$   and $a,b:[\alpha,\beta]\to\mathbb{R}$ are continuous functions.
Consider a polynomial $F(z)=\sum\limits_{j=0}^{n}\delta_j z^j$,  $\delta_j\in\mathbb{R}$, $j=0,\ldots,n$.
Then,
$  AB=BF(A)$ if and only if \
$
  \supp (a-F(a(\gamma))) \,\cap  \supp b=\emptyset.
$
\end{prop:DjinjaTumwSilv1}

\begin{proof}
We have
\begin{gather*}
  (ABx)(t)=A(Bx)(t)=a(t)(Bx)(t)=a(t)b(t)x(\gamma) \\
  (A^mx)(t)=a(t)^{m} x(t),\quad m=0,\ldots,n\\
  (BA^m)x)(t)=(BA^mx)(t)=B(A^mx)(t)=\left\{\begin{array}{rc}
   b(t)x(\gamma), &  m=0    \\
   b(t)  a(\gamma)^{m} x(\gamma),& \quad m=1,\ldots, n
  \end{array}\right. \\
  (BF(A)x)(t)= \delta_0 b(t)x(\gamma)+\sum\limits_{j=1}^{n} \delta_j b(t) a(\gamma)^{j} x(\gamma)=b(t)F(a(\gamma))x(\gamma).
\end{gather*}
Then $(ABx)(t)=(BF(A)x)(t)$ if and only if for all $x\in C[\alpha,\beta]$ and for all $t\in [\alpha,\beta]$,
\begin{equation*}
  a(t)b(t)x(\gamma)=b(t)F(a(\gamma))x(\gamma).
\end{equation*}
 This is equivalent to the equation $ a(\cdot)b(\cdot)=b(\cdot)F(a(\gamma))$,
 which is equivalent to
 \begin{equation*}
   \supp\, (a-F(a(\gamma))) \,\cap  \supp\, b=\emptyset,
 \end{equation*}
 by Lemma \ref{LemmaContinuousFunctionsSupportMeasureZero}. \qed
\end{proof}

\begin{ex:DjinjaTumwSilv1}{\rm
Let $A:C[\alpha,\beta]\to C[\alpha,\beta]$,  $B:C[\alpha,\beta]\to C[\alpha,\beta]$ be defined by
\begin{equation*}
  (Ax)(t)=a(t)x(t) ,\quad (Bx)(t)= b(t) x(\gamma),
\end{equation*}
where $\alpha,\beta\in\mathbb{R}$, $\alpha<\beta$, $\gamma\in [\alpha,\beta]$   and $a,b:[\alpha,\beta]\to\mathbb{R}$ are continuous functions.
Then, for a real constant $\zeta$ and a positive integer $d$,
\begin{equation*}
 AB=\zeta BA^d
 \end{equation*}
 if and only if
$
  \supp (a-\zeta a(\gamma)^d) \,\cap  \supp b=\emptyset.
$
This follows by Proposition  \ref{PropAmultOpBinitialValueOpContFunct}.
}\end{ex:DjinjaTumwSilv1}

\begin{ex:DjinjaTumwSilv1}{\rm
Let $A=A_\nu:C[0,2]\to C[0,2]$,  $B:C[0,2]\to C[0,2]$ be defined by
\begin{equation*}
  (Ax)(t)=A_{\nu}x(t)=a(t)x(t) ,\quad (Bx)(t)= b(t) x(1/2),
\end{equation*}
where
$
a(t)=I_{[0,1]}(t)\sin(\pi t)+\nu$, $ b(t)=I_{[1,2]}(t)\sin(\pi t)$,
$\nu\in\mathbb{R}$, and $I_{[\alpha_1,\beta_1]}(t)$ is the indicator function of the interval $[\alpha_1,\beta_1]$. Consider a  polynomial $F(z)=\sum\limits_{j=0}^{n}\delta_j z^j$ with  $\delta_i\in\mathbb{R}$ for $i=0,\ldots,n$. If
$ \nu=F(1+\nu),
$
then, by Proposition \ref{PropAmultOpBinitialValueOpContFunct}, these operators satisfy  $AB=BF(A)$.
Note that the condition $\nu=F(1+\nu)$ is equivalent to $\tilde{F}(1+\nu)=1+\nu$, where $\tilde{F}(z)= F(z)+1$, that is to
$1+\nu$ being a fixed point
of $\tilde{F}(z)$, or equivalently to $1+\nu$ being the root of
$\tilde{F}(z)-z=0$, which in terms of $F$ is the same as $1+\nu$ being a root of $F(z)-z+1=0$.
If this equation has roots in $\mathbb{R}$, then for such roots the corresponding operators $A$ and $B$ satisfy $AB=BF(A)$. If $F(z)=z-1$, then $1+\nu$ is a root of  $F(z)-z+1=0$ for any real number $\nu$, since the equation becomes the equality $0=0$. Thus, for any $\nu\in \mathbb{R}$, the operators $A$ and $B$ satisfy
$AB=B(A-I)$ which is equivalent to $AB-BA=-B$ and to $BA-AB=B$. This can be also checked by direct computation. In fact, the operator $A_{\nu}$ corresponding to the parameter $\nu$ can be expressed as  $A_{\nu}=A_{0}+\nu I$. For $\nu=0$, $a$ and $b$ have non-overlapping supports and hence $(A_{0} Bx)(t)=x(\frac{1}{2})a(t)b(t)=0$,
and $(BA_{0} x)(t)=b(t)(A_{0}x)(\frac{1}{2})=b(t) a(\frac{1}{2})x(\frac{1}{2})=b(t)x(\frac{1}{2})=B$, and thus $BA-AB=B$, and hence, for any $\nu\in \mathbb{R}$,
$BA_{\nu}-A_{\nu}B=B(A_{0}+\nu I)-(A_{0}+\nu I)B=BA_{0}+\nu B-A_{0}B-\nu B=BA_{0}-A_{0}B=B$.
}\end{ex:DjinjaTumwSilv1}

\begin{cor:DjinjaTumwSilv1}
Let $A:C[\alpha,\beta]\to C[\alpha,\beta]$,  $B:C[\alpha,\beta]\to C[\alpha,\beta]$ be defined by
\begin{equation*}
  (Ax)(t)=a(t)x(t) ,\quad (Bx)(t)= b(t) x(\gamma),
\end{equation*}
where $\alpha,\beta\in\mathbb{R}$, $\alpha<\beta$, $\gamma\in [\alpha,\beta]$   and $a,b:[\alpha,\beta]\to\mathbb{R}$ are continuous functions.
Consider a polynomial $F(z)=\sum\limits_{j=0}^{n}\delta_j z^j$,  $\delta_j\in\mathbb{R}$, $j=0,\ldots,n$.
If $a(\gamma)=0$ then $  AB=BF(A)$ if and only if
$ \supp (a-\delta_0) \,\cap  \supp b=\emptyset.$
Furthermore, if $\delta_0\not =0$ then $AB=BF(A)$ yields $b(\gamma)=0$.
\end{cor:DjinjaTumwSilv1}

\begin{proof}
  This follows by Proposition \ref{PropAmultOpBinitialValueOpContFunct}.
   \qed
\end{proof}


\acknowledgement{
This work was supported by the Swedish International Development Cooperation
Agency (Sida) bilateral program with  Mozambique. Domingos Djinja is grateful to the research environment Mathematics and Applied Mathematics (MAM), Division of Mathematics and Physics at M{\"a}lardalen University, School of Education, Culture and Communication,
M\"alardalens University for excellent environment for research in Mathematics.
Partial support from Swedish Royal Academy of Sciences is also gratefully acknowledged.}

\end{document}